\documentclass[a4paper,12pt]{article}
\usepackage[frenchb,english]{babel}
\usepackage{amssymb}
\usepackage{amsmath}

\usepackage{hyperref}
\usepackage{pb-diagram}
\usepackage{graphicx}
\newtheorem{thm}{Theorem}[section]
\newtheorem{cor}[thm]{Corollary}
\newtheorem{lemme}[thm]{Lemma}
\newtheorem{rema}[thm]{Remark}

\newtheorem{pro}[thm]{Proposition}

\newcommand{\iid}{\mathrm{Id}\,}
\newcommand{\pfend}{$\hfill\square$}
\newcommand{\Ss}{\mathbb{S}}
\newcommand{\HH}{\mathbb{H}}
\newcommand{\beqt}{\begin{equation}}  \newcommand{\eeqt}{\end{equation}}
\newcommand{\beQ}{\begin{eqnarray*}} \newcommand{\eeQ}{\end{eqnarray*}}

\newcommand{\ddet}{\mathrm{det}}
\newcommand{\R}{\mathbb{R}}

\begin{document}
\title{{\bf Skew Killing spinors}}
\author{{\small\bf{Georges Habib}}\footnote{Lebanese University, Faculty of Sciences II, Department of Mathematics, P.O. Box 90656 Fanar-Matn, Lebanon, E-mail: ghabib@ul.edu.lb} , {\small\bf{ Julien Roth}}\footnote{LAMA, Universit\'e Paris-Est Marne-la-Vall\'ee, Cit\'e Descartes, Champs-sur-Marne, F-77454 Marne-la-Vall\'ee, Cedex 2, France, E-mail: julien.roth@univ-mlv.fr}}\date{}
\maketitle

\begin{abstract} 
\noindent In this paper, we study the existence of a skew Killing spinor (see the def\mbox{}inition below) on $2$ and $3$-dimensional Riemannian spin  manifolds. We establish the integrability conditions  and prove that these spinor f\mbox{}ields correspond to twistor spinors in the two dimensional case while, up to a conformal change of the metric, they correspond to parallel spinors in the  three dimensional case.
\end{abstract} 

\section{Introduction} 
\setcounter{equation}{0}
Classifying spin manifolds $(M^n,g)$ carrying particular spinor f\mbox{}ields has interested many authors, since it is related to several geometric constructions (see \cite{Bo,Gi} for results in this topic). A $\lambda$-Killing spinor is a smooth section $\psi$ of the spinor bundle $\Sigma M$ satisfying, for all $X\in \Gamma(TM)$, the differential equation $\nabla_X\psi=\lambda X\cdot\psi$ where $\lambda$ is a complex number, ``$\cdot$''and $\nabla$ denote respectively the Clifford multiplication and the spinorial Levi-Civita connection on $\Sigma M$. 
The existence of such a spinor imposes very rigid conditions on the geometry of the manifold. In particular, $M$ is Einstein with constant scalar curvature equal to $4\lambda^2n(n-1)$. As a consequence, $\lambda$ should be either real or purely imaginary. In \cite{B1}, C. B\"ar gave a geometric description for simply connected manifolds with real Killing spinors (i.e. those correspond to $\lambda$ real and different from zero). He showed that in fact, there exists a one-to-one correspondence between Killing spinors on $\Sigma M$ and parallel spinors on the cone of $M$, obtained by the warped product of the manifold $M$ with the half real line. Therefore his classif\mbox{}ication is based on Wang's one for complete and simply connected manifolds carrying parallel spinors \cite{Wa}. Manifolds with imaginary Killing spinors (i.e. for $\lambda$ purely imaginary) have been treated by H. Baum \cite{BFGK} in a different way by studying the foliation def\mbox{}ined by the length function of the spinor. The problem of classifying manifolds with spinor f\mbox{}ields $\psi$ satisfying for all $X\in \Gamma(TM)$, the equation 
\begin{equation}\label{eq:01}
\nabla_X\psi=A(X)\cdot\psi
\end{equation} 
where $A$ is an endomorphism tensor f\mbox{}ield on $TM$ (they are $\lambda$-Killing spinors for $A=\lambda {\rm Id}$) was open until now. First, we point out that the symmetric part of $A$ is the energy-momentum tensor of $\psi$ studied in \cite{H3}. In \cite{Fr1}, Th. Friedrich proved that, in the two dimensional case, if there exists a solution of Equation $\eqref{eq:01}$ with a symmetric tensor $A$, then $A$ satisf\mbox{}ies the Codazzi-Mainardi and Gauss equations. Consequently, the surface $M$ is locally immersed into the $3$-Euclidean space with Weingarten tensor equal to $-2A.$ Conversely, the restriction of any parallel spinor on $\mathbb{R}^3$ to an oriented hypersurface is a solution of Equation \eqref{eq:01}. 
We notice that fundamental topological informations can be read off from those spinor representations \cite{KS}. This statement was later generalized in \cite{B2} for all dimensions (see also \cite{Mo}). In fact, the authors showed that, for a symmetric and Codazzi tensor $A$ (i.e. $(\nabla A)(X,Y)=(\nabla A)(Y,X)$ for all $X,Y\in \Gamma(TM)$), there is a one-to-one correspondence between solutions of \eqref{eq:01} and parallel spinors on the generalised warped product endowed with a deformation of the metric in terms of $A.$\\ 

\noindent In this paper, we aim at studying the solutions of Equation \eqref{eq:01} for low dimensions, where $A$ is a non-trivial endomorphism tensor f\mbox{}ield. In the case of a skew-symmetric tensor, we call these {\it skew Killing spinors}. \\
First, we consider the two dimensional case. We show that any surface carrying such a spinor can be locally immersed into a 3-dimensional round sphere with Weingarten tensor given by the symmetric part of $A.$ As a corollary, we get some rigidity results for surfaces with particular $A$. Moreover, we prove that these spinors correspond, in particular cases, to twistor spinors \cite{BFGK}, i.e. spinor f\mbox{}ields $\varphi$ satisfying for all $X\in \Gamma(TM)$, the equation $\nabla_X\varphi=-\frac{1}{2}X\cdot D\varphi$ where $D$ denotes the Dirac operator on $M$. This allows us to describe them on the round sphere. \\
We also consider the case of imaginary tensor $iA$ (see Equation \eqref{eq:02}) 
which can be obtained from local isometric immersions in the anti-de Sitter space $\HH^{2,1}$.\\
In Section $3,$ we discuss the three dimensional case. Indeed, we show if $M$ admits a skew Killing spinor, then the tangent bundle of $M$ splits into two integrable bundles, corresponding to the kernel of $A$ and its orthogonal. Moreover, the integrability conditions (see Equations \eqref{eq:con}) imply that $M$ is conformally f\mbox{}lat and using a suitable conformal change of the metric on the universal cover $\widetilde M$, we show that they are in correspondence with parallel spinors. 
To illustrate our results, we treat several important examples.

\section{The 2-dimensional case}
\setcounter{equation}{0}
\noindent This section is devoted to the study of solutions of Equation \eqref{eq:01} on $2$-dimensional manifolds. We also consider the case where $A$ is an imaginary endomorphism tensor field.
We begin with the real case.
\subsection{Case of real tensor}
We start from the following fact: If $(M^2,g)$ is a surface of $\Ss^3,$ endowed with its standard metric of curvature one, the restriction of a real Killing spinor on $M$ is a spinor field $\psi$ satisfying the following equation 
\beqt\label{eqkilling2}
\nabla_X\psi=-\frac{1}{2}S(X)\cdot\psi-\frac{1}{2}X\cdot\omega\cdot\psi,
\eeqt
where $S$ is the Weingarten tensor of the surface and $\omega$ stands the real volume form given by $\omega=e_1\cdot e_2$. Moreover, B. Morel \cite{Mo} showed that the existence of such a spinor on a surface $(M,g)$ is equivalent to the existence of a local isometric immersion of $M$ into $\Ss^3$.\\
On the other hand, if we denote by $J$ the complex structure of $M$ given by the rotation of angle $\frac{\pi}{2}$ on $TM$, it is obvious that $X\cdot\omega=-J(X)$. Hence, Equation \eqref{eqkilling2} reduces to 
$$
\nabla_X\psi=-\frac{1}{2}S(X)\cdot\psi+\frac{1}{2}J(X)\cdot\psi.
$$
In other terms, the spinor field $\psi$ is a solution of \eqref{eq:01} with $A=-\frac{1}{2}(S-J)$.
We now prove the following result.

\begin{thm}\label{thm1}
Let $(M^2,g)$ be a connected and oriented Riemannian surface and let $A$ be a field of endomorphism on $TM$. We set $A=S+T$, with $S$ symmetric and $T$ skew-symmetric. If there exists on $M$ a spinor field satisfying Equation \eqref{eq:01}
with a Codazzi tensor $S$. Then $T=bJ$ with $b$ a real constant and there exists a local isometric immersion of $M$ into the sphere $\Ss^3(4b^2)$ with Weingarten tensor $-2S$.
\end{thm}
\begin{rema}\label{equDirac}
The existence of such a spinor is equivalent to the existence of a spinor of constant norm solution of the Dirac equation
$$D\psi=H\psi+2b\omega\cdot\psi,$$
where $H$ denotes the mean curvature of the surface.
\end{rema}

\noindent{\bf Proof.} A direct computation of the spinorial curvature of $\psi$ gives
$$\mathcal{R}(X,Y)\psi=d^{\nabla}A(X,Y)+\big(A(Y)\cdot A(X)-A(X)\cdot A(Y)\big)\cdot\psi.$$
By the Ricci identity \cite[p.156]{Bo}, one gets
$$-\frac{1}{2}R_{1212}\omega\cdot\psi=\mathcal{R}(e_1,e_2)\psi=d^{\nabla}A(e_1,e_2)\cdot\psi-2\ddet(A)\omega\cdot\psi.$$
Thus, since $\ddet(A)=\ddet(S)+\ddet(T)$, we have
$$\Big(R_{1212}-4\ddet(S)-4\ddet(T)\Big)\omega\cdot\psi=-2d^{\nabla}A(e_1,e_2)\cdot\psi.$$
Using the classical argument in \cite{Fr1} and \cite{Mo}, we obtain
\beqt\label{GaussCodazzi}
R_{1212}-4\ddet(S)-4\ddet(T)=0\quad\text{and}\quad  d^{\nabla}A(e_1,e_2)=0.
\eeqt
The tensor $A$ is then a Codazzi tensor which leads to a Codazzi $T$, by assumption on $S$. On the other hand, the tensor $T$ is skew-symmetric, so it is of the form $bJ$ with $b$ a real-valued function. Thus, we get
\beQ
0=d^{\nabla}T(e_1,e_2)&=&\nabla_{e_1}(T(e_2))-\nabla_{e_2}(T(e_1))\\
&=&-e_1(b)e_1-e_2(b)e_2=-\nabla b.
\eeQ
Since $M$ is connected, the function $b$ is constant. Then, Equation \eqref{GaussCodazzi} becomes
the Gauss and Codazzi equations for an isometric immersion into the sphere of curvature $4b^2$ with Weingarten tensor given by $-2S$. 
\pfend\\ \\
Actually, the hypothesis on $S$ to be a Codazzi tensor is somehow restrictive on the surface. For this, we need to introduce two spinor fields in order to conclude.
\begin{pro}
Let $(M^2,g)$ be a connected and oriented Riemannian surface and let $A$ be a field of endomorphism on $TM$. We set $A=S+T$, with $S$ symmetric and $T$ skew-symmetric and $B=S-T$. If there exists on $M$ two spinor fields $\varphi$ and $\psi$ satisfying
$$\nabla_X\varphi=A(X)\cdot\varphi\quad\text{and}\quad \nabla_X\psi=B(X)\cdot\psi,$$
then $T=bJ$ with $b$ constant and there exists a local isometric immersion of $(M^2,g)$ into the sphere $\Ss^3(4b^2)$ and Weingarten tensor $-2S$.
\end{pro}
{\bf Proof.} From the computation of the proof of Theorem \ref{thm1}, we easily see that the tensors $A$ and $B$ are both Codazzi. Since $2S=A+B$, we deduce that $S$ is Codazzi and we conclude using Theorem \ref{thm1}. 
\pfend\\ \\
Now when the Weingarten tensor is particular, we may deduce some rigidity results on the surface. In the following, we treat the case where $S=a{\rm Id}$ for some real-valued function $a$. 
\begin{cor}\label{cor:SKS}
Let $(M^2,g)$ be a connected and oriented Riemannian surface. If $M$ carries a spinor field $\psi$ such that
\begin{equation}\label{twi}
\nabla_X\psi=a X\cdot\psi+bJ(X)\cdot\psi,
\end{equation}
with $a,b$ are real-valued functions. If $a$ is a real constant, then $b$ is constant and $M$ is isometric to the sphere $\Ss^2(4a^2+4b^2)$.
\end{cor}
{\bf Proof.} From the fact that $S=a\iid$ using Theorem \ref{thm1}, we know that $M$ is immersed into $\Ss^3(b^2)$ and is totally umbilical. Hence, $M$ is a geodesic sphere of  $\Ss^3(b^2)$ with mean curvature $a$ and so $M$ is isometric to a $2$-dimensional sphere of curvature $4a^2+4b^2$. 
\pfend
\begin{cor}
 Let $(M^2,g)$ be an oriented Riemannian surface. If $M$ carries a spinor fields $\psi$ such that
 $$
\nabla_X\psi=a X\cdot\psi+aJ(X)\cdot\psi,
$$
where $a$ is a non-trivial real-valued function. Then $a$ is constant and $M$ is isometric to the sphere $\Ss^2(2a^2)$.
\end{cor}
{\bf Proof.} Since the tensor $A=a\iid+aJ$ is Codazzi, we compute
\beQ 0=d^{\nabla}A(e_1,e_2)&=&\nabla_{e_1}(ae_2-ae_1)-\nabla_{e_2}(ae_1+ae_2)\\
&=&e_1(a)e_2-e_1(a)e_1-e_2(a)e_1-e_2(a)e_2.
\eeQ
Therefore, we get that $e_1(a)=e_2(a)=0$ and so $a$ is constant. From Theorem \ref{thm1}, we deduce the existence of an isometric immersion of  $(M^2,g)$ into $\Ss^3(4a^2)$ with Weingarten tensor $-2a\iid$. Hence, $(M^2,g)$ is isometric to a sphere of curvature $2a^2$.
\pfend\\

\noindent On the other hand, using the remark above, one may easily check that spinor fields satisfying Equation \eqref{twi} are actually twistor spinors of constant norm, since $\Re(\nabla_X\psi,\psi)=\frac{1}{2}X(|\psi|^2)$ for all $X\in \Gamma(TM)$. Now conversely, we will prove that any twistor spinor of constant norm can satisfy such an equation. We have the following proposition:
\begin{pro} \label{pro:twis} Let $(M^2,g)$ be an oriented Riemannian manifold carrying a twistor spinor $\psi$ of norm $1$. Then there exists two real valued functions $a$ and $b$ on $M$ such that Equation \eqref{twi} holds for $\psi$.
\end{pro} 
{\bf Proof.} Since the spinor f\mbox{}ield $\psi$ is of norm $1$ and the spinor bundle of $M$ is of real rank $4$, the covariant derivative of $\psi$ can be expressed with respect to the orthonormal frame $\{\psi,e_1\cdot\psi,e_2\cdot\psi, e_1\cdot e_2\cdot\psi\}$ for all $X\in \Gamma(TM)$ as 
\begin{eqnarray}
\nabla_X\psi&=&\Re(\nabla_X\psi,\psi)\psi+\Re(\nabla_X\psi,e_1\cdot\psi)e_1\cdot\psi+\Re(\nabla_X\psi,e_2\cdot\psi)e_2\cdot\psi\nonumber\\
&&+\Re(\nabla_X\psi,e_1\cdot e_2\cdot\psi)e_1\cdot e_2\cdot\psi.
\label{eq:2}
\end{eqnarray}
The f\mbox{}irst term in the r.h.s. of \eqref{eq:2} clearly vanishes from the fact that the norm of $\psi$ is constant. On the other hand, the fourth term also vanishes. Indeed, using that $Z\cdot\nabla_Z\psi$ is independent of the choice of any vector f\mbox{}ield $Z$ of norm $1$ \cite{BFGK}, one has for $X=e_1$ (analogously for $X=e_2$)  
\begin{eqnarray*}
\Re(\nabla_{e_1}\psi,e_1\cdot e_2\cdot \psi)&=&-\Re(e_1\cdot\nabla_{e_1}\psi,e_2\cdot\psi)\\ 
&=&-\Re(e_2\cdot\nabla_{e_2}\psi,e_2\cdot\psi)
=-\frac{1}{2}e_2(|\psi|^2)=0.
\end{eqnarray*}
Thus in order to conclude, it is sufficient to compute the components of the tensor $A(X,Y)=\Re(\nabla_X\psi,Y\cdot\psi)$ with respect to the frame $\{e_1,e_2\}.$ Indeed, it is a direct fact that they are equal to $A(e_i,e_i)=-\frac{1}{2}\Re(D\psi,\psi)$ for $i=1,2$ and that 
\begin{eqnarray*} 
A(e_1,e_2)&=&\Re(\nabla_{e_1}\psi,e_2\cdot\psi) 
=-\Re(e_1\cdot e_2\cdot\nabla_{e_2}\psi,e_2\cdot\psi)\\ 
&=&-\Re(\nabla_{e_2}\psi,e_1\cdot\psi)=-A(e_2,e_1),
\end{eqnarray*}
from which follows the proof of the proposition. \hfill$\square$\\\\ 

\noindent In order to justify the existence of solutions of Equation \eqref{twi} with constant $a$ on the round sphere, it is sufficient by \ref{pro:twis} to find a twistor spinor $\psi$ of constant norm satisfying $\Re(D\psi,\psi)={\rm cte}$. For this, let $\varphi$ and $\theta$ be orthogonal Killing spinors of Killing constants equal respectively to $\frac{1}{2}$ and $-\frac{1}{2}$. Therefore, the spinor f\mbox{}ield $\psi:=\varphi+\theta$ is a twistor spinor of constant norm $|\varphi|^2+|\theta|^2$, since it is the sum of two Killing spinors of different constants \cite{BFGK}. Furthermore, we have $D\psi=-\varphi+\theta$ and thus the real product $\Re(D\psi,\psi)=|\theta|^2-|\varphi|^2$ is constant. 

\subsection{Case of imaginary tensor}
In this section, we study the existence of spinor fields which are solutions of the following equation
\beqt\label{eq:02}
\nabla_X\psi=iA(X)\cdot\psi,
\eeqt
where $A$ is an arbitrary tensor field. Here again, such spinors appear as restriction of parallel or Killing spinors but of Lorentzian space forms.
In \cite{LR}, the second author and M.A. Lawn proved an analogue to Friedrich's result for isometric immersions of Riemannian surfaces into Lorentzian space forms. First, the restriction of a parallel spinor of the Minkowski space $\R^{2,1}$ to a Riemannian surface is a solution of \eqref{eq:02}  
where $-2A$ is the Weingarten tensor of the immersion. Conversely, the existence of such a spinor on $M$ with a Codazzi tensor $A,$ implies the existence of an isometric  immersion of $(M,g)$ into $\R^{2,1}$. Note that the assumption on $A$ to be Codazzi is a necessary fact. 
Similarly, the restriction of an imaginary Killing spinor of the anti-de Sitter space $\HH^{2,1}$ to a Riemannian surface gives a spinor $\psi$ verifying
\beQ
\nabla_X\psi&=&-\frac{i}{2}S(X)\cdot\psi-\frac{i}{2}X\cdot\omega\cdot\psi\\
&=&-\frac{i}{2}S(X)\cdot\psi+\frac{i}{2}J(X)\cdot\psi.
\eeQ
We show the following result:
\begin{thm}\label{thm2}
Let $(M^2,g)$ be a connected and oriented Riemannian surface and $A$ a field of endomorphism on $TM$. We set $A=S+T$, with $S$ symmetric and $T$ skew-symmetric. We assume that $S$ is Codazzi. Then, if there exists on $M$ two orthogonal and nowhere vanishing spinor fields $\varphi$ and $\psi$ 
$$\nabla_X\varphi=iA(X)\cdot\varphi\quad\text{and}\quad \nabla_X\psi=iA(X)\cdot\psi,$$
then $T=bJ$ with $b$ constant and there exists an isometric immersion of $(M^2,g)$ into $\HH^{2,1}(4b^2)$ with Weingarten tensor $-2S$.
\end{thm}
{\bf Proof.} The scheme of proof is similar to Theorem \ref{thm1}. The computation of the spinorial curvature of $\varphi$ and $\psi$ gives
$$\Big(\underbrace{R_{1212}+4\ddet(S)+4\ddet(T)}_{G}\Big)\omega\cdot\varphi=-2i\underbrace{d^{\nabla}A(e_1,e_2)}_{C}\cdot\varphi,$$
and
$$\Big(\underbrace{R_{1212}+4\ddet(S)+4\ddet(T)}_{G}\Big)\omega\cdot\psi=-2i\underbrace{d^{\nabla}A(e_1,e_2)}_{C}\cdot\psi.$$
Hence, we deduce
\beQ
<C\cdot\varphi,\psi>&=&<-iC\cdot\varphi,-i\psi>\\
&=&\frac{1}{2}<G\omega\cdot\varphi,-i\psi>=\frac{1}{2}<\varphi,iG\omega\cdot\psi>\\
&=&<\varphi,C\cdot\psi>=-<C\cdot\varphi,\psi>.
\eeQ
Then $C\cdot\varphi$ and $\psi$ orthogonal. But, by assumption, $\varphi$ and $\psi$ are orthogonal. Since the spinor bundle is of complex rank $2$, the spinor field $C\cdot\varphi$ and $\varphi$ are colinear, and there exists a complex-valued function $f$ so that $C\cdot\varphi=f\varphi$. Taking the real scalar product by $\varphi$, we see that $f$ can only have imaginary values. Hence  $f=ih$ with $h$ real-valued. So, we have $G\omega\cdot\varphi=2h\varphi$. We take again the real scalar product by $\varphi$, which yields to
$$2h|\varphi|^2=G<\omega\cdot\varphi,\varphi>=0.$$
Since $\varphi$ never vanishes, we deduce that $h$ vanishes everywhere. Therefore, $G=0$ et $C=0$. The fact that $C=0$ implies that $A$ is Codazzi since $S$ is Codazzi, then $T$ too. Hence, $T=bJ$ with $b$ constant. The  equation $G=0$ becomes
$$R_{1212}+4\ddet(S)+4b^2=0.$$
Thus, the Gauss and Codazzi equations for an isometric immersion into $\HH^{2,1}(4b^2)$ are fulfilled. 
\hfill$\square$ \\ \\
As in the real case, we have a correspondence between solutions of \eqref{eq:02} and twistor spinors. Indeed, we prove the following result.
\begin{pro} Let $(M^2,g)$ be an oriented Riemannian manifold carrying a twistor spinor $\psi$ so that $\Re(\psi,\overline{\psi})=1$. Then there exists two real valued functions $a$ and $b$ on $M$ such that 
$$\nabla_X\psi=iaX\cdot\psi+ibJ(X)\cdot\psi.
$$
\end{pro} 
{\bf Proof.}  The proof is similar to the one in Proposition \ref{pro:twis} using the fact that $\left\{\frac{\overline{\psi}}{|\psi|},e_1\cdot\frac{\overline{\psi}}{|\psi|},e_2\cdot\frac{\overline{\psi}}{|\psi|}, e_1\cdot e_2\cdot\frac{\overline{\psi}}{|\psi|}\right\}$ is an orthonormal frame.
\hfill$\square$
\noindent\section{The 3-dimensional case}
\setcounter{equation}{0}
In this section, we will study the existence of skew Killing spinors in the $3$-dimensional case. As mentioned in the introduction, we will prove that the tangent bundle of $M$ splits into two integrable subbundles. Moreover, we will show that skew Killing spinors correspond to parallel spinors, up to a conformal change of the metric. 
\begin{pro} \label{pro:decom}Let $(M^3,g)$ be an oriented Riemannian manifold. Assume that $M$ carries a skew Killing spinor $\psi$. Then we have an orthogonal decomposition  
$$TM=L\oplus Q $$ 
where $L:={\rm Ker} A,$ and $Q:=L^\perp$ is an integrable normal bundle. 
\end{pro}
{\bf Proof.} The decomposition follows directly from the fact that $0$ is a simple eigenvalue for any skew-symmetric matrix. We shall now prove that $Q$ is integrable. For this, let us choose a direct orthonormal frame $\{\xi,e_1,e_2\}$ of $\Gamma(TM)$ such that $\xi\in \Gamma(L)$ and $\{e_1,e_2\}$ in $\Gamma(Q)$. With respect to this frame, the matrix $A$ can be written as $\begin{pmatrix}0&0&0\\0&0&-b\\0&b&0\end{pmatrix},$ where $b$ is a real-valued function on $M.$ In this case, Equation \eqref{eq:01} is then reduced to
\begin{equation}
\nabla_\xi\psi=0, \,\, \nabla_{e_1}\psi=be_2\cdot\psi, \,\, \nabla_{e_2}\psi=-be_1\cdot\psi.
\label{eq:Kill}
\end{equation}  
As in the two dimensional case, we aim to compute the spinorial curvature $R(\cdot,\cdot)\psi$ in the local frame $\{\xi,e_1,e_2\}$. For this, we denote in the sequel by $\kappa:=\nabla_\xi\xi,$ and by $h:=\nabla\xi$ the endomorphism f\mbox{}ield in $\Gamma(Q)$. The different components of the curvature are then equal to
\begin{eqnarray*}
R(\xi,e_1)\psi&=&\nabla_\xi\nabla_{e_1}\psi-\nabla_{e_1}\nabla_{\xi}\psi-\nabla_{[\xi,e_1]}\psi\\ 
&=&\nabla_\xi(be_2\cdot\psi)-g([\xi,e_1],e_1)\nabla_{e_1}\psi-g([\xi,e_1],e_2)\nabla_{e_2}\psi\\
&=&-bg(\kappa,e_2)\xi\cdot\psi-bg(h(e_1),e_2)e_1\cdot\psi+(\xi(b)+bg(h(e_1),e_1)) e_2\cdot\psi.
\end{eqnarray*} 
Similarly, we compute 
\begin{eqnarray*}
R(\xi,e_2)\psi
&=&bg(\kappa,e_1)\xi\cdot\psi-(\xi(b)+bg(h(e_2),e_2))e_1\cdot\psi+bg(h(e_2),e_1) e_2\cdot\psi.
\end{eqnarray*} 
Also we have 
\begin{eqnarray*}
R(e_1,e_2)\psi
&=&(-2b^2+b{\rm Trace}(h))\xi\cdot\psi-e_1(b)e_1\cdot\psi-e_2(b) e_2\cdot\psi.
\end{eqnarray*} 
Using the Ricci identity and the fact that the complex volume form $\omega_3=-\xi\cdot e_1\cdot e_2$ acts as the identity on $\Sigma M$, one gets for the Ricci curvatures 
\begin{eqnarray}\label{eq:Ric1} 
-\frac{1}{2}{\rm Ric}(e_1)\cdot \psi&=&e_2\cdot R(e_1,e_2)\psi+\xi\cdot R(e_1,\xi)\psi\nonumber\\
&=& (e_2(b)-bg(\kappa,e_2))\psi+e_1(b)\xi\cdot\psi\nonumber\\
&&+(-2b^2+b{\rm Trace}(h)+bg(h(e_1),e_1)+\xi(b))e_1\cdot\psi\nonumber\\&&+bg(h(e_1),e_2)e_2\cdot\psi. 
\end{eqnarray}
Further with an analogous computation, we f\mbox{}ind 
\begin{eqnarray}
-\frac{1}{2}{\rm Ric}(e_2)\cdot \psi
&=&(-e_1(b)+bg(\kappa,e_1))\psi+e_2(b)\xi\cdot\psi+bg(h(e_2),e_1)e_1\cdot\psi\nonumber\\
&&+(-2b^2+b{\rm Trace}(h)+bg(h(e_2),e_2)+\xi(b))e_2\cdot\psi,
\label{eq:Ric2}
\end{eqnarray} 
and, 
\begin{eqnarray} 
-\frac{1}{2}{\rm Ric}(\xi)\cdot \psi
&=&b(g(h(e_1),e_2)-g(h(e_2),e_1))\psi+(2\xi(b)+b{\rm Trace}(h))\xi\cdot\psi\nonumber\\
&&+bg(\kappa,e_1)e_1\cdot\psi+bg(\kappa,e_2)e_2\cdot\psi.
\label{eq:Ric3}
\end{eqnarray}
The Hermitian product of Equations \eqref{eq:Ric1}, \eqref{eq:Ric2} and \eqref{eq:Ric3} by $\psi$ gives after identifying the real parts, the following conditions 
\begin{equation}
\left\{\begin{array}{ll}
g(h(e_1),e_2)=g(h(e_2),e_1)\\\\
e_1(b)=bg(\kappa,e_1)\\\\
e_2(b)=bg(\kappa,e_2) 
\end{array}\right.
\label{eq:con}
\end{equation} 
which imply the integrability of the normal bundle, since we have $g([e_1,e_2],\xi)=g(h(e_2),e_1)-g(h(e_1),e_2)=0.$ 
\hfill$\square$\\\\
We point out that the equations above for the Ricci curvatures with \eqref{eq:con} are in fact of global type. Recall f\mbox{}irst that the frame $\{\xi,e_1,e_2\}$ is direct, so if one interchanges $\xi$ onto $\xi':=-\xi,$ since it also belongs to the kernel of $A$, one should replace the function $b$ by $-b$. Hence the equations above remain the same. However, by a suitable modif\mbox{}ication of the Levi-Civita connection to a flat connection, one can prove that these equations are suff\mbox{}icient to def\mbox{}ine a skew Killing spinor. We have 
\begin{pro} \label{pro:cond}Let $(M^3,g)$ be a simply connected Riemannian manifold and $b$ a real-valued function on $M$ and let $\{\xi,e_1,e_2\}$ be a direct orthonormal frame such that Equations \eqref{eq:Ric1}, \eqref{eq:Ric2}, \eqref{eq:Ric3} and \eqref{eq:con} are satisf\mbox{}ied. Then there exists a skew Killing spinor. 
\end{pro}  
In the following, we will see that the existence of a skew Killing spinor gives rise to a particular geometry of the manifold.  We recall that a Riemannian manifold $(M^n,g)$ is said to be conformally f\mbox{}lat if for each point $x\in M$ there exists a neighbourhood $V$of $x$ and a smooth function $u$ such that $(V,e^{2u})$ is f\mbox{}lat. This is equivalent to saying, in the three dimensional case, that the Schouten tensor $S:=\frac{{\rm Scal}}{4}g-{\rm Ric}$ has the property $(\nabla_X S)Y=(\nabla_Y S)X$ for all vector f\mbox{}ields $X,Y\in \Gamma(TM)$ \cite{BFGK}. Here ${\rm Scal}$ denotes the scalar curvature of $M$. In the following proposition, we will see that the orthogonal decomposition of $TM$ established in Proposition \ref{pro:decom} shows that $M$ is conformally f\mbox{}lat.    
\begin{pro} \label{pro:con} Let $(M^3,g)$ be an oriented Riemannian manifold carrying a skew Killing spinor $\psi$, then $M$ is conformally flat.
\end{pro}
{\bf Proof.} Let $\{\xi,e_1,e_2\}$ be an orthonormal frame considered as in Proposition \ref{pro:decom}. We aim to show that the above property for $S$ is satisf\mbox{}ied for the vector f\mbox{}ields $X=e_1$ and $Y=e_2$. The same computations can be done for the other vector f\mbox{}ields. Indeed, the Dirac operator associated with $\psi$ is equal to 
\begin{equation}
D\psi=\xi\cdot\nabla_\xi\psi+e_1\cdot\nabla_{e_1}\psi+e_2\cdot\nabla_{e_2}\psi=2be_1\cdot e_2\cdot\psi=2b\xi\cdot\psi.
\label{eq:6}
\end{equation} 
First, it is straightforward to see by Equations \eqref{eq:Ric1}, \eqref{eq:Ric2} and \eqref{eq:Ric3}, that the scalar curvature of $M$ is equal to $8(b^2-\xi(b)-b{\rm Trace}(h)).$ On the other hand, by differentiating \eqref{eq:6}, it follows that
\begin{eqnarray*}
\nabla_{e_1}D\psi&=&2e_1(b)\xi\cdot\psi+2bh(e_1)\cdot\psi+2b^2\xi\cdot e_2\cdot\psi\nonumber\\ 
&=&2e_1(b)\xi\cdot\psi+(2bg(h(e_1),e_1)-2b^2)e_1\cdot\psi+2bg(h(e_1),e_2)e_2\cdot\psi\nonumber\\
&=&-{\rm Ric}(e_1)\cdot\psi-2(b{\rm Trace}(h)+\xi(b)-b^2)e_1\cdot\psi\nonumber\\
&=&-{\rm Ric}(e_1)\cdot\psi+\frac{1}{4}{\rm Scal}\,e_1\cdot\psi=S(e_1)\cdot\psi,
\end{eqnarray*}
also that $\nabla_{e_2}D\psi=S(e_2)\cdot\psi$ and $\nabla_{\xi}D\psi=S(\xi)\cdot\psi-2b^2\xi\cdot\psi.$
Differentiating again, one gets 
$$\nabla_{e_1}\nabla_{e_2}D\psi=\nabla_{e_1}(S(e_2))\cdot\psi+bS(e_2)\cdot e_2\cdot\psi,$$
and
$$\nabla_{e_2}\nabla_{e_1}D\psi=\nabla_{e_2}(S(e_1))\cdot\psi-bS(e_1)\cdot e_1\cdot\psi,$$
and $\nabla_{[e_1,e_2]}D\psi=S([e_1,e_2])\cdot\psi$. The curvature can then be written as  
$$R(e_1,e_2)D\psi=\{(\nabla_{e_1}S)(e_2)-(\nabla_{e_2}S)(e_1)\}\cdot\psi+bS(e_2)\cdot e_2\cdot\psi+bS(e_1)\cdot e_1\cdot\psi.$$
Using again \eqref{eq:6} and the fact that  
\begin{equation*}
R(e_1,e_2)(\xi \cdot\psi)=R(e_1,e_2,e_2,e_1)\psi-\xi\cdot R(e_1,e_2)\psi,
\end{equation*} 
which can be proved by a direct computation, we get by linearity of the curvature that 
\begin{eqnarray*}
2b(R(e_1,e_2,e_2,e_1)\psi-\xi\cdot R(e_1,e_2)\psi)&=&\{(\nabla_{e_1}S)(e_2)-(\nabla_{e_2}S)(e_1)\}\cdot\psi\\&&-b(S(e_2,e_2)+S(e_1,e_1))\psi\\&&-bS(e_2,\xi) e_1\cdot\psi+ bS(e_1,\xi)e_2\cdot\psi.
\end{eqnarray*} 
Inserting the formula $R(e_1,e_2,e_2,e_1)=\frac{1}{2}(\sum_{i=1}^2{\rm Ric}(e_i,e_i)-{\rm Ric}(\xi,\xi)),$ into the equation above
we deduce with the help of Proposition \ref{pro:decom} that
\begin{eqnarray*} 
-b{\rm Ric}(\xi,\xi)\psi+2b(-2b^2+b{\rm Trace}(h))\psi&=&\{(\nabla_{e_1}S)(e_2)-(\nabla_{e_2}S)(e_1)\}\cdot\psi\\&&-\frac{b}{2}{\rm Scal}\psi,
\end{eqnarray*}
which gives $(\nabla_{e_1}S)(e_2)=(\nabla_{e_2}S)(e_1)$. This completes the proof.
\hfill$\square$\\\\ 
Now, we will state the following lemma which allows us to determine the conformal factor of the metric. 
\begin{lemme} Under the same hypotheses as in Proposition \ref{pro:con}, the vector f\mbox{}ield $\tau:=b\xi$ is globally def\mbox{}ined and is a closed $1$-form.
\end{lemme}
{\bf Proof.} The f\mbox{}irst point is clear since as we said before, when one chooses the vector f\mbox{}ield $\xi'=-\xi$ then Equations \eqref{eq:Kill} still hold, with respect to the frame $\{\xi',e_2,e_1\},$ with the function $b'=-b$. In both cases $\tau$ is globally def\mbox{}ined. Now we compute the exterior derivative of $\tau$
\begin{eqnarray*}
d\tau&=&\xi\wedge \nabla_\xi(b\xi)+e_1\wedge\nabla_{e_1}(b\xi)+e_2\wedge\nabla_{e_2}(b\xi)\\
&=&b\xi\wedge\kappa+e_1(b)e_1\wedge\xi+b e_1\wedge h(e_1)+e_2(b)e_2\wedge\xi+b e_2\wedge h(e_2)\\
&=&bg(\kappa,e_1)\xi\wedge e_1+bg(\kappa,e_2)\xi\wedge e_2+e_1(b)e_1\wedge\xi+bg(h(e_1),e_2)e_1\wedge e_2\\&&+e_2(b)e_2\wedge\xi+bg(h(e_2),e_1) e_2\wedge e_1=0, 
\end{eqnarray*}
from which we deduce the lemma with the help of Equations \eqref{eq:con}. $\hfill\square$\\\\ 
Since the vector f\mbox{}ield $\tau$ is a closed $1$-form there exists on the universal cover $\widetilde M$ a real-valued function $u$ such that $du=-2\tau$. We will in fact see that a skew Killing spinor on $M$ induces a parallel spinor on  $(\widetilde M,\overline{g}:=e^{2u}g)$. First, let us review some standard facts on conformal metrics. We recall that a given spin structure on $(M,g)$ induces a spin structure on $(M,\overline{g})$ for any conformal change of the metric $\overline{g}$. Moreover, there exists an isometry (with respect to the Hermitian product) between the corresponding spinor bundles $\Sigma_g M$ and $\Sigma_{\overline g}M$. Together with the corresponding isometry of the tangent bundle, given by $X\rightarrow \overline{X}:=e^{-u}X,$ the Clifford multiplications $``\cdot"$ and $``\,\overline\cdot\,"$ are then related for every $X\in \Gamma(TM)$ and $\psi \in \Gamma(\Sigma_g M)$ by the formula $\overline{X}\,\overline{\cdot}\,{\overline\psi}=\overline{X\cdot\psi}$. Here ${\overline\psi}$ denotes the spinor f\mbox{}ield associated with $\psi$ with respect to the isometry \cite{Bo}. Now we have the following proposition 
\begin{pro} \label{pro:paral}Let $(M^3,g)$ be an oriented Riemannian manifold carrying a skew Killing spinor $\psi$, then the universal cover $(\widetilde M, \overline{g})$ carries a parallel spinor $\overline\psi$.    
\end{pro}
{\bf Proof.} The proof is based on applying to $\psi$ the relation between the Levi-Civita connections $\nabla$ and $\overline\nabla$ on the spinor bundles $\Sigma_g M$ and $\Sigma_{\overline g}M$. We write for all $X\in \Gamma(TM)$ \cite{Bo} 
\begin{equation}
\overline\nabla_X\overline\psi=\overline{\nabla_X\psi}-\frac{1}{2}\overline{X\cdot du\cdot\psi}-\frac{1}{2}X(u)\overline\psi.
\label{eq:confo}
\end{equation}
Using Equations \eqref{eq:Kill} and the fact that $X(u)=-2bg(\xi,X)$ yield in the local frame $\{\xi,e_1,e_2\},$  
$$\overline\nabla_{\xi}\overline\psi=b\overline{\xi\cdot\xi\cdot\psi}+b\overline{\psi}=0$$
and also 
$$\overline\nabla_{e_1}\overline\psi=b\overline{e_2\cdot\psi}+b\overline{e_1\cdot \xi\cdot\psi}=0.$$
The same computation can be done for $e_2$ and we f\mbox{}inish the proof of the proposition. \hfill$\square$\\\\
Conversely, we shall see that any parallel spinor on a three Riemannian manifold def\mbox{}ines a skew Killing spinor for any conformal change of the metric. We have 
\begin{pro} \label{pro:confor} Let $(M^3,g)$ be a Riemannian manifold with a parallel spinor $\psi$, then for any conformal change $\overline{g}=e^{2u}g$ where $u:M\rightarrow \mathbb{R}$ is a real-valued function, the manifold $(M,\overline{g})$ carries a skew Killing spinor $\overline\psi$. 
\end{pro}
{\bf Proof.} As in the previous proposition, we apply the formula \eqref{eq:confo} to the spinor f\mbox{}ield $\psi$ and we get
$${\overline \nabla}_X\overline \psi=-\frac{1}{2}\overline{X\cdot du\cdot\psi}-\frac{1}{2}X(u)\overline\psi,$$
for all $X\in \Gamma(TM)$. Now for any direct orthonormal frame $\{e_1,e_2,e_3\}$ on $(M,g)$, it follows
\begin{eqnarray*} 
{\overline \nabla}_{\overline{e}_1}\overline \psi&=&-\frac{1}{2}e^{-u}\overline{e_1\cdot du\cdot \psi}-\frac{1}{2}\overline{e}_1(u)\overline\psi\\
&=&-\frac{1}{2}e^{-u}\overline{e_1\cdot (e_1(u)e_1+e_2(u)e_2+e_3(u)e_3)\cdot \psi}-\frac{1}{2}\overline{e}_1(u)\overline\psi\\
&=& -\frac{1}{2}\overline{e}_2(u)\overline{e}_3\overline{\cdot}\overline{\psi}+\frac{1}{2}\overline{e}_3(u)\overline{e}_2\overline{\cdot}\overline{\psi}.
\end{eqnarray*} 
The same computations can be done for $\overline{e}_2$ and $\overline{e}_3$ and thus
$${\overline \nabla}_{\overline{e}_2}\overline \psi=\frac{1}{2}\overline{e}_1(u)\overline{e}_3\overline{\cdot}\overline{\psi}-\frac{1}{2}\overline{e}_3(u)\overline{e}_1\overline{\cdot}\overline{\psi},$$ 
also 
$${\overline \nabla}_{\overline{e}_3}\overline \psi=-\frac{1}{2}\overline{e}_1(u)\overline{e}_2\overline{\cdot}\overline{\psi}+\frac{1}{2}\overline{e}_2(u)\overline{e}_1\overline{\cdot}\overline{\psi}.$$ 
Hence the equations above can be written in a homogeneous way for all $X\in \Gamma(TM)$ as ${\overline \nabla}_{X}\overline \psi=A(X)\overline{\cdot} \overline{\psi}$ where $A$ is the endomorphism f\mbox{}ield on $\Gamma(TM)$ given by $A(X)=-\frac{1}{2}*_{\overline g}(d_{\overline g}u\wedge X)$, which is skew-symmetric, since 
$${\overline g}(AX,Y)v_{\overline g}=-\frac{1}{2}d_{\overline g}u\wedge X\wedge Y=\frac{1}{2}d_{\overline g}u\wedge Y\wedge X=-{\overline g}(AY,X)v_{\overline g},$$ 
here $v_{\overline g}$ denotes the volume form associated with the metric ${\overline g}$. This f\mbox{}inishes the proof of the proposition. \hspace{85mm}$\square$\\\\
From the above propositions, we conclude that skew Killing spinors are in a bijective correspondance with parallel spinors. The existence of a parallel spinor on a Riemannian manifold implies the vanishing of the Ricci curvature which is equivalent, in the three dimensional case, to the f\mbox{}latness of the metric. Simply connected f\mbox{}lat $3$-manifolds are globally immersed into the Euclidean space $\mathbb{R}^3$, with the local isometry property. However, the completeness of the $\overline{g}$ in Proposition \ref{pro:paral} is not assured as we shall see later in Example $2$. \\\\
On the other hand, for any skew Killing spinor $\psi$ def\mbox{}ined on a Riemannian spin manifold $(M,g)$, one may associate the vector f\mbox{}ield $\zeta$ def\mbox{}ined for all $X\in \Gamma(TM)$ by $g(\zeta,X)=i(X\cdot\psi,\psi)$. A straightforward computation for the covariant derivative of $\zeta$ leads to 
$$\nabla_X\zeta=2g(\zeta,\tau)X-2g(X,\zeta)\tau,$$
where we recall that $\tau=b\xi$ is a closed one form. An immediate consequence is that the norm of $\zeta$ is constant and that $d\zeta=-2\zeta\wedge\tau$. Hence the couple $(\zeta,-2\tau)$ is the Pfaff form on $M$ def\mbox{}ining a transversally aff\mbox{}ine orientable foliation of codimension $1$ given by the equation $\zeta=0$ \cite{Se}. These foliations are locally def\mbox{}ined by submersions over the Euclidean space $\mathbb{R}$ where the germs of the transverse coordinates belong to the real aff\mbox{}ine group. The universal cover of $M$ is then a submersion over $\mathbb{R}$. \\

\noindent We f\mbox{}inish the paper by treating some examples admitting skew Killing spinors. These are global products of manifolds of one dimension with manifolds of two dimensions, endowed with some particular metrics. Our technical computations are based on applying Proposition \ref{pro:cond} that we leave the details to the reader.\\ 

\noindent{\bf Example 1.} Let consider the manifold $(M,g)=(\mathbb{R}^2\times \mathbb{R}, \theta^2g_{{\rm stan}}\oplus f^2dz^2)$ where $\theta:\mathbb{R}\rightarrow \mathbb{R}$ is a function of $z$ and $f:\mathbb{R}^3\rightarrow \mathbb{R}$ is a function of $x,y,z$. The set $\{\xi:=\frac{1}{f}\partial z, e_1:=\frac{1}{\theta}\partial x, e_2:=\frac{1}{\theta}\partial y\}$ is direct and orthonormal frame with respect to the metric $g$. In the sequel, we will use the notation ${\dot f}_x$ for the partial derivative of $f$ with respect to the variable $x$. It is easy to show using Koszul's formula that the tensor $h|_{\xi^\perp}=\nabla\xi$ is given by the matrix   
$\begin{pmatrix}\frac{\dot \theta(z)}{\theta f}&0\\0&\frac{\dot \theta(z)}{\theta f}\end{pmatrix}$ with respect to the frame $\{e_1,e_2\}$ and the term $\kappa=\nabla_\xi\xi$ has coordinates  $(-\frac{{\dot f}_x}{\theta f},-\frac{{\dot f}_y}{\theta f})$. Therefore, it is suff\mbox{}icient to f\mbox{}ind a function $b$ such that Equations \eqref{eq:Ric1}, \eqref{eq:Ric2}, \eqref{eq:Ric3} and \eqref{eq:con} are satisf\mbox{}ied. Now a long computation for the Ricci curvatures yields after comparing that


$$b=\frac{\dot \theta(z)}{2\theta f} \quad\text{and}\quad f(x,y,z)=c(z)x+d(z)y+e(z),$$
for real functions $c,d,e.$\\

\noindent Consider the direct frame $\{\xi=\frac{1}{\theta}\partial x,\,\, e_1=\frac{1}{f}\partial y,\,\, e_2=\frac{1}{f}\partial z\}$ and we shall see that there exists also a skew Killing spinor with respect to this frame. Similar computations as above for the tensor $h$ and $\kappa$ give after comparing the Ricci curvatures give 
that $b=\frac{B(x)}{\theta(z)}$ with $B$ is a function of $x$ and $8$ other relations which lead us to distinguish two cases:\\\\ 
{\bf Case where $\theta(z)=1$}: In this case, we f\mbox{}ind that $B$ is constant and thus $b$ is constant. Moreover, 
$$f(x,y,z)=\frac{c(z)}{2b}e^{2bx}+d(z){\rm cos}(2by)+e(z){\rm sin}(2by),$$ 
where $d,e$ and $c$ are functions of $z$. The function $f$ descends to the product $\mathbb{T}^2\times \mathbb{R}$ if and only if $c(z)=0$. If moreover the functions $d(z),e(z)$ are periodic, the function $f$ descends to the torus $\mathbb{T}^3$.\\\\ 
{\bf Case where $\dot f_y=0$}: We have that $b=\frac{\dot f_x}{2\theta f}$ with the equations: 
$$-\ddot \theta(z) f+\dot \theta(z)\dot f_z=0,\,\, \dot \theta(z)^2=\ddot f_x f-\dot f_x^2,\,\, -\ddot f_{xz}+\dot f_x\dot f_z=0.$$  
The solution of the above differential system gives that $\theta(z)=cz+d$ with $c,d\in \mathbb{R}$ and $f(x,y,z)=Ae^{\sqrt{H}x}+Be^{-\sqrt{H}x}$ with $H>0, A,B\in \mathbb{R}$ and $4ABH=c^2$ and $c\neq 0$.\\    

\noindent{\bf Example 2.} Consider the manifold $(M,g)=(\mathbb{R}\times \mathbb{S}^2, dt^2\oplus f^2g_{{\rm stan}}),$ where $f:\mathbb{R}\rightarrow \mathbb{R}$ is a function of $t$ and denote by $\xi$ the vector f\mbox{}ield $\partial t$ and $\{e_1,e_2\}$ an orthonormal frame on $(\mathbb{S}^2, f^2g_{{\rm stan}})$.
First we recall that $\nabla_{\xi}\xi=0$ and $h|_{\xi^\perp}=\frac{\dot{f}}{f}g_t$ \cite[ch.7]{O}. For the Ricci curvatures, we have for all vector fields $X,Y\in \Gamma(T\mathbb{S}^2)$ that 
${\rm Ric}(X,\xi)=0,\,\, {\rm Ric}(\xi,\xi)=-2\frac{\ddot{f}}{f}$ and   
\begin{equation}
{\rm Ric}(X,Y)= \frac{1-\dot{f}^2-f\ddot{f}}{f^2}g_t(X,Y),
\label{eq:Ric}
\end{equation}
From Equations \eqref{eq:con}, we deduce f\mbox{}irst that $b$ is a function of $t$. Moreover, by \eqref{eq:Ric3} and \eqref{eq:Ric1}, we get the relation $(\dot{f}-2bf)^2=1$.\\\\
For the particular solution $f(t)=1$ and $b(t)=\frac{1}{2}$, the metric $\overline{g}=e^{-2t}(dt^2\oplus g_{{\rm stan}})$ carries a parallel spinor by Propositon \ref{pro:paral}, but it is not a complete metric, as we mentionned before. Indeed, we f\mbox{}ix a point $x$ on the sphere and we consider the sequence $u_n=(n,x)$ for $n\in \mathbb{N}$. The sequence $(u_n)_n$ is clearly divergent but it is in fact a Cauchy sequence, since for $p>q$ the distance 
$d(u_p,u_q)=\frac{1}{2}(e^{-2q}-e^{-2p})$
tends to zero when $p,q\rightarrow \infty$.\\\\ 
\noindent{\bf Example 3.} Let $\{e_1,e_2,e_3\}$ be three vector f\mbox{}ields on $\mathbb{R}^3$ given by the Lie brackets 
$$[e_1,e_2]=0,\,\, [e_1,e_3]=e_1,\,\, [e_2,e_3]=e_2.$$
We consider a metric $g$ on $\mathbb{R}^3$ such that the frame $\{e_1,e_2,e_3\}$ is orthonormal. The Christofel symbols $\Gamma_{ij}^k=g(\nabla_{e_i}e_j,e_k)$ are then equal to 
$$\Gamma_{11}^3=\Gamma_{22}^3=-\Gamma_{13}^1=-\Gamma_{23}^2=-1,$$
and all others $\Gamma_{ij}^k$ are equal to zero. 
Taking a constant spinor f\mbox{}ield $\psi$, we can easily prove by using the def\mbox{}inition of the covariant derivative that 
$$\nabla_{e_3}\psi=0,\,\,\nabla_{e_1}\psi=\frac{1}{2}e_2\cdot\psi,\,\, \nabla_{e_2}\psi=-\frac{1}{2}e_1\cdot\psi.$$ 

\end{document}